\documentclass[12pt]{article}
\usepackage{color}
\definecolor{darkblue}{rgb}{0.00,0.25,0.50}
\usepackage[colorlinks,filecolor=blue,citecolor=darkblue]{hyperref}

\usepackage{amsfonts,amssymb,amsmath,amsthm}
\usepackage{url}
\usepackage{enumerate}
\usepackage[ukrainian, russian, english]{babel}
\usepackage[cp1251]{inputenc}

\usepackage{cite}
\begin{document}
\selectlanguage{ukrainian}
\thispagestyle{empty}

\title{}

\begin{center}
\textbf{\Large Наближення сумами Зигмунда класів згорток періодичних функцій  в інтегральних метриках}
\end{center}
\vskip0.5cm
\begin{center}
 У.\,З.~Грабова\\ \emph{\small Східноєвропейський національний університет імені Лесі Українки, Луцьк}
\end{center}
\vskip0.5cm

\begin{abstract}
We obtain estimates exact in order for deviations  of Zygmund sums in metrics of spaces  $L_{q}$, $1<q<\infty$, on classes    of   $2\pi$-periodic functions, that admit the representation in the form of convolution of functions that belong to unit ball of the space $L_{1}$ with fixed kernel $\Psi_{\beta}$.   We show that at certain values of the parameters that define the class $L^{\psi}_{\beta,1}$ and method of approximation, Zygmund sums  provide the order of best approximation of given classes by trigonometric polynomials in metric $L_{q}$.

\vskip0.6cm

Получены точные по порядку оценки отклонений сумм Зигмунда в метриках пространств  $L_{q}$, $1<q<\infty$, на классах    $2\pi$-периодических функций, которые допускают изображение в виде свертки функций,  принадлежащих единичному шару пространства  $L_{1}$ с фиксированным ядром $\Psi_{\beta}$. Показано, что при определенных значениях параметров, что определяют класс   $L^{\psi}_{\beta,1}$ и метод приближения, суммы Зигмунда обеспечивают порядок наилучшего приближения указанных классов тригонометрическими полиномами в метрике  $L_{q}$.
\end{abstract}

\vskip1cm


Нехай $L_{p}$, $1\leq p\leq\infty$, --- простір  $2\pi$-періодичних сумовних функцій $\varphi$   зі
скінченною нормою   $\| \varphi\|_{p}$, де при $p\in [1,\infty)$
$
\|\varphi\|_{p}:=\Big(\int\limits_{0}^{2\pi}|\varphi(t)|^{p}dt\Big)^{\frac{1}{p}}$,
а при $p=\infty$ ${\|\varphi\|_{\infty}:=\mathop{\rm{ess}\sup}\limits_{t}|\varphi(t)|}$.

Нехай, далі,  $L^{\psi}_{\beta,p}$, \ $1\leq p<\infty$,  --- клас
$2\pi$-періодичних функцій $f(x)$, що майже для всіх
$x\in\mathbb{R}$ представляються згортками

\begin{equation}\label{zgo}
f(x)=\frac{a_{0}}{2}+(\Psi_{\beta}\ast\varphi)(x):=\frac{a_{0}}{2}+
\frac{1}{\pi}\int\limits_{-\pi}^{\pi}\Psi_{\beta}(x-t)\varphi(t)dt,
\ a_{0}\in\mathbb{R}, \ \varphi\perp 1,
\end{equation}
де  $\|\varphi \|_{p}\leq 1$,  а $\Psi_{\beta}(t)$ --- сумовне на
$[0,2\pi)$ ядро, ряд Фур'є якого має вигляд
$$
\sum\limits_{k=1}^{\infty}\psi(k)\cos
\big(kt-\frac{\beta\pi}{2}\big), \ \psi(k)>0,\ \beta\in
    \mathbb{R}.
$$
 Функцію $\varphi$ в зображенні (\ref{zgo}) називають
$(\psi,\beta)$-похідною функції $f$ і позначають через
$f^{\psi}_{\beta}$. Поняття $(\psi,\beta)$-похідної належить О.І.
Степанцю (див., наприклад, \cite[с.~132]{S1}).

Якщо $\psi(k)=k^{-r}$, $r>0$, $\beta\in  \mathbb{R}$, то класи $L^{\psi}_{\beta,p}$ є  класами
Вейля-Надя і  позначаються через $W^{r}_{\beta,p}$.

Не зменшуючи загальності можна вважати, що послідовність $\psi(k)$, яка визначає класи згорток, є слідом на множині натуральних чисел $\mathbb{N }$ деякої
неперервної функції $\psi(t)$ неперервного аргументу $t\geq1$.

Позначимо через $\Theta_{\varrho}$, $1\leq \varrho<\infty$,
 множину монотонно незростаючих
функцій $\psi(t)$, $t\geq1$, для яких існує стала $\alpha>\frac{1}{\varrho}$
така, що функція $t^{\alpha}\psi(t)$ майже спадає,  тобто для будь-яких $t_{1}>t_{2}\geq 1$ виконується нерівність $t^{\alpha}_{1}\psi(t_{1})\leq Kt^{\alpha}_{2}\psi(t_{2})$, в якій
 $K$ --- деяка додатна стала. Прикладами функцій $\psi$, що задовольняють умову  $\psi\in
\Theta_{\varrho}$ є, зокрема, функції вигляду
${\psi_{1}(t)=\frac{1}{t^{r}}}$, \ $r>\frac{1}{\varrho}$;
$\psi_{2}(t)=\frac{\ln^{\alpha}(t+c)}{t^{r}}$, \ $r>\frac{1}{\varrho}$,
\ $\alpha>0$, \ $c>e^{\frac{2\alpha}{r-\frac{1}{\varrho}}}-1$;
$\psi_{3}(t)=\frac{1}{t^{r}\ln^{\alpha}(t+c)}$, \
$r>\frac{1}{\varrho}$, \ $\alpha>0$, \ $c>0$; \  $\psi_{4}(t)=\frac{\ln
\ln^{\alpha}(t+c)}{t^{r}}$, \ $r>\frac{1}{\varrho}$, \ $\alpha>0$, \
$c>e^{\frac{2\alpha}{r-\frac{1}{\varrho}}}-1$.

 У даній роботі будемо вимагати, щоб  $\psi\in\Theta_{q'}$, $\frac{1}{q}+\frac{1}{q'}=1$, $1<q<\infty$.
 Умова $\psi\in \Theta_{q'}$, як неважко переконатись,
гарантує справедливість включення $\Psi_{\beta}\in L_{q}$ (див., наприклад, \cite[с.~657]{Bari}), а отже, і вкладення $L^{\psi}_{\beta,1}\subset
L_{q}$, $1<q<\infty$ (див., наприклад, нерівність (1.5.28) із роботи \cite[с.~43]{Korn}).

Сумами Зигмунда
 функції  $f$  із  $L_{1}$ називають тригонометричні
поліноми вигляду
\begin{equation}\label{sz}
Z_{n-1}^{s}(f;t)=\frac{a_{0}}{2}+\sum\limits_{k=1}^{n-1}\bigg(1-\Big(\frac{k}{n}\Big)^{s}\bigg)(a_{k}\cos
kt+b_{k}\sin kt), s>0,
\end{equation}
де $a_{k}$ і $b_{k}$ --- коефіцієнти Фур'є функції $f$.
При $s=1$ поліноми $Z_{n-1}^{s}$ є відомими сумами Фейєра, які
позначаються через $\sigma_{n-1}(f;t)$.

В даній роботі досліджуються величини
$$
{\cal E}\big(L^{\psi}_{\beta,p};
Z_{n-1}^{s}\big)_{{q}}=\mathop{\sup}\limits_{f\in
L^{\psi}_{\beta,p}}\|f(\cdot)-Z^{s}_{n-1}(f;\cdot)\|_{{q}},
$$
з метою одержання для них точних порядкових оцінок, при  $\psi\in \Theta_{q'}$, $\beta\in \mathbb{R}$, $s>0$, $1<q<\infty$ і $p=1$.

 В роботі також розглядатимуться найкращі наближення класів  $L^{\psi}_{\beta,p}$, $p\geq1$, $\beta\in  \mathbb{R}$
тригонометричними поліномами $t_{n-1}$ порядку $n-1$, тобто величини вигляду
$$
E_{n}\big(L^{\psi}_{\beta,p}
\big)_{q}=\sup\limits_{f\in
L^{\psi}_{\beta,p}}\inf\limits_{t_{n-1}}\|f(\cdot)-t_{n-1}(\cdot)\|_{q}, \ 1\leq p, \  q\leq\infty,
$$
за умови, що $L^{\psi}_{\beta,p}\subset
L_{q}$.
Далі буде вказано область допустимих значень параметрів $\psi$, $s$ і $q$ при яких порядки спадання величин ${\cal E}\big(L^{\psi}_{\beta,1};
Z_{n-1}^{s}\big)_{{q}}$ і $E_{n}\big(L^{\psi}_{\beta,1}
\big)_{q}$ збігаються.

Вивченню апроксимативних властивостей сум Зигмунда на різних функціональних класах присвячена значна кількість робіт (з детальною бібліографією з цього напряму можна ознайомитись в \cite{SG}). Зокрема, у  \cite{SG}  знайдено порядкові оцінки величин ${\cal E}\big(L^{\psi}_{\beta,p};
Z_{n-1}^{s}\big)_{{\infty}}$ при $\psi\in \Theta_{p}$, $1<p<\infty$ і $\beta\in \mathbb{R}$.

Порядки спадання  величин   $E_{n}\big(L^{\psi}_{\beta,p}
\big)_{q}$
 досліджувались багатьма авторами (див. \cite{S2,UMG} та ін.). У \cite{UMG} встановлено точні порядкові рівності величин $E_{n}\big(L^{\psi}_{\beta,1}
\big)_{q}$, \linebreak$1<q\leq\infty$, за умови, що $\psi\in \Theta_{q'}\cap B$, $\frac{1}{q}+\frac{1}{q'}=1$ і $\beta\in \mathbb{R}$, де $B$ --- множина монотонно незростаючих функцій $\psi(t)$, для кожної з яких можна вказати додатну сталу $K$ таку, що
$
 \frac{\psi(t)}{\psi(2t)}\leq K,\   \ t\geq1.
$

Дана робота є продовженням досліджень, розпочатих у \cite{SG} та \cite{UMG}.
Для формулювання основних результатів введемо наступні означення.

Будемо казати, що додатна функція $g(t)$, задана на $[1,\infty)$,
належить до множини $A^{+}$ ($g\in A^{+}$), якщо існує $\varepsilon>0$ таке, що
$g(t)t^{-\varepsilon}$ зростає на $[1,\infty)$.  Якщо
існує ${\varepsilon>0}$ таке, що $g(t)t^{\varepsilon}$ спадає на
$[1,\infty)$, то будемо казати, що $g$ належить до множини
$A^{-}$ ($g\in A^{-}$). Через $\mathcal{Z}$ позначимо множину неперервних слабо
коливних (в сенсі Зигмунда) функцій, тобто  додатних функцій
$g(t)$, визначених на $[\frac{1}{\pi},\infty)$, таких, що при довільному
$\delta>0$  для достатньо великих $t$ \   $g(t)t^{\delta}$ зростає, а $g(t)t^{-\delta}$
спадає.

 Далі під записом $A(n)=O(B(n))$ розумітимемо, що
існує стала $K>0$, така, що виконується нерівність
$A(n)\leq K(B(n))$, для всіх $n\in \mathbb{N}$. Запис
$A(n)\asymp B(n)$ означає, що $A(n)=O(B(n))$ і одночасно
$B(n)=~O(A(n))$.

\textbf{Теорема 1.}
\emph{Нехай $1< q<\infty$, $s>0$,  $n\in \mathbb{N}$,
$g_{s,q'}(t):=\psi(t)t^{s+\frac{1}{q'}}$, $\frac{1}{q}+\frac{1}{q'}=1$ і $\beta\in \mathbb{R}$. Тоді}

\emph{1. Якщо $\psi\in \Theta_{q'}$ і $g_{s,q'}\in A^{+}$, то}
\begin{equation}\label{tp5}
{\cal E}\big(L^{\psi}_{\beta,1};
Z_{n-1}^{s}\big)_{q}=O\big(\psi(n)n^{1-\frac{1}{q}}\big).
\end{equation}

 \emph{2. Якщо   $g_{s,q'}\in \mathcal{Z}$,
   то}

\begin{equation}\label{tp4}
{\cal E}\big(L^{\psi}_{\beta,1};
Z_{n-1}^{s}\big)_{q}=O\bigg(\frac{1}{n^{s}}
\Big(\int\limits_{1}^{n}\frac{\big(g_{s,q'}(t)\big)^{q}}{t}dt\Big)^{\frac{1}{q}}\bigg).
\end{equation}

\emph{3. Якщо  $g_{s,q'}\in A^{-}$,
  то}
\begin{equation}\label{tp1}
{\cal E}\big(L^{\psi}_{\beta,1};
Z_{n-1}^{s}\big)_{q}=O\big(n^{-s}\big).
\end{equation}

\textbf{\emph{Доведення.}}  З рівностей (\ref{zgo}) та (\ref{sz}) випливає, що для будь-якої
$f\in L^{\psi}_{\beta,1}$ майже для всіх $x\in \mathbb{R}$ має місце
рівність
\begin{equation}\label{predst1p}
f(x)-Z_{n-1}^{s}(f;x)=\frac{1}{\pi}\int\limits_{-\pi}^{\pi}\Big(\frac{1}{n^{s}}\sum\limits_{k=1}^{n-1}
\psi(k)k^{s}\cos\Big(k(x-t)-\frac{\beta\pi}{2}\Big)\!+\!\Psi_{\beta,n
    }(x-t)\Big)\varphi(t)dt,
\end{equation}
де $\|\varphi\|_{1}\leq1$, $\Psi_{\beta,n
    }(\tau)=\sum\limits_{k=n}^{\infty}\psi(k)\cos\Big(k\tau-\frac{\beta\pi}{2}\Big)$,
\ \   $n\in \mathbb{N}$.

Із  (\ref{predst1p}), застосовуючи нерівність (1.5.28) з роботи \cite [с.~43]{Korn} та нерівність трикутника, отримуємо
$$
{\cal E}\big(L^{\psi}_{\beta,1};
Z_{n-1}^{s}\big)_{q}
\leq\frac{1}{\pi n^{s}}\bigg\|
\sum\limits_{k=1}^{n-1}\psi(k)k^{s}\cos\Big(k(\cdot)-\frac{\beta\pi}{2}\Big)\bigg\|_{q}+
\frac{1}{\pi}\big\|\Psi_{\beta,n}(\cdot)\big\|_{q}.
$$

Тоді, скориставшись  формулами (18), (21), (24) та (42) з роботи \cite{SG},  одержимо оцінки  \eqref{tp5}--\eqref{tp1}.
Теорему 1 доведено.

\textbf{Теорема 2.} \emph{Нехай $1< q<\infty$, $s>0$,  $n\in \mathbb{N}$,
$g_{s,q'}(t):=\psi(t)t^{s+\frac{1}{q'}}$, $\frac{1}{q}+\frac{1}{q'}=1$ і $\beta\in \mathbb{R}$.
  Тоді}

\emph{1. Якщо $\psi\in \Theta_{q'}$, $g_{s,q'}\in A^{+}$ і функція $1/\psi(t)$ опукла вгору або донизу на $[1,\infty)$, то}
\begin{equation}\label{tp52}
{\cal E}\big(L^{\psi}_{\beta,1};
Z_{n-1}^{s}\big)_{q}\asymp E_{n}\big(L^{\psi}_{\beta,1}
\big)_{q}\asymp\psi(n)n^{1-\frac{1}{q}}.
\end{equation}

\emph{2. Якщо   $g_{s,q'}\in \mathcal{Z}$, то}

\begin{equation}\label{tp2zgama}
{\cal E}\big(L^{\psi}_{\beta,1};
Z_{n-1}^{s}\big)_{q}\asymp\frac{1}{n^{s}}
\Big(\int\limits_{1}^{n}\frac{\big(g_{s,q'}(t)\big)^{q}}{t}dt\Big)^{\frac{1}{q}}.
\end{equation}

\emph{3. Якщо   $g_{s,q'}\in A^{-}$, то}
\begin{equation}\label{tp12p}
{\cal E}\big(L^{\psi}_{\beta,1};
Z_{n-1}^{s}\big)_{q}\asymp n^{-s}.
\end{equation}

\textbf{\emph{Доведення.}} Оцінки зверху для величин ${\cal E}\big(L^{\psi}_{\beta,1};
Z_{n-1}^{s}\big)_{q}$ в (\ref{tp52}) -- (\ref{tp12p}) випливають із співвідношень (\ref{tp5}) -- (\ref{tp1}) відповідно.

Розглянемо випадок $g_{s,q'}\in A^{+}$, $\psi\in \Theta_{q'}$.
У роботі \cite[с.~1195]{UMG}, при  $\psi\in
\Theta_{q'}\cap B$, $1<q\leq\infty$, $\frac{1}{q}+\frac{1}{q'}=1$, за умови, що функція $1/\psi(t)$ опукла вгору або донизу на $[1,\infty)$, встановлено наступну оцінку     знизу для величини
${E}_{n}(L^{\psi}_{\beta,1})_{q}$:
\begin{equation}\label{poroznnnp}
{E}_{n}(L^{\psi}_{\beta,1})_{q}\geq
K_{\psi,q}\psi(n)n^{1-\frac{1}{q}}, \ n\in \mathbb{N}, \ \beta\in
\mathbb{R},
\end{equation}
де $K_{\psi,q}$ --- додатна величина, що може залежати лише від
$\psi$ та $q$. Оскільки з включення $g_{s,q'}\in A^{+}$, $s>0$, $1< q\leq\infty$  випливає, що $\psi\in B$,  то з оцінок  (\ref{tp5}), (\ref{poroznnnp}) та нерівності \ \
$
{E}_{n}(L^{\psi}_{\beta,1})_{q}\leq{\cal
E}\big(L^{\psi}_{\beta,1};Z_{n-1}^{s}\big)_{q}, \ n\in
\mathbb{N}
$
\ \  випливає порядкова рівність (\ref{tp52}).

Нехай, далі, $g_{s,q'}\in \mathcal{Z}$. Для встановлення оцінки знизу величини
${\cal E}\big(L^{\psi}_{\beta,1};
Z_{n-1}^{s}\big)_{q}$, $1<q<\infty$ розглянемо функцію
$$
\varphi_{\alpha}(t)=\varphi_{\alpha}(n,t)=\alpha \Big(V_{n}(t)-\frac{1}{2}\Big), \ \alpha>0, \ n\in
\mathbb{N},
$$
де $V_{m}(t)$ --- ядра методу Валле Пуссена,
$$
V_{m}(t):=\frac{1}{m}\sum\limits_{k=m}^{2m-1}D_{k}(t)=D_{m}(t)+2\sum\limits_{k=m+1}^{2m-1}\Big(1-\frac{k}{2m}\Big)\cos kt, \ m\in \mathbb{N},
$$
$D_{k}(t)$ --- ядра Діріхле, \ \
$
D_{k}(t):=\frac{1}{2}+\sum\limits_{\nu=1}^{k}\cos \nu t=\frac{\sin\big(k+\frac{1}{2}\big)t}{2\sin\frac{t}{2}}, \ k\in  \mathbb{N}.
$

Оскільки  (див., наприклад, \cite[с.~1192]{UMG}) $|V_{m}(t)|<A_{1}m$, $0\leq t\leq\pi$, \ $|V_{m}(t)|\leq\frac{A_{2}}{mt^{2}}$, $0<t\leq\pi$, де $A_{i}$ --- абсолютні сталі, то нескладно переконатися, що
\begin{equation}\label{ochfi}
\|\varphi_{\alpha}\|_{1}\leq \alpha\big(\pi+\|V_{n}\|_{1}\big)\leq\alpha A_{3}.
\end{equation}
При $\alpha=\alpha_{0}=A_{3}^{-1}$ з (\ref{ochfi}) одержуємо  нерівність
$\|\varphi_{\alpha_{0}}\|_{1}\leq1$.

 Із \cite[c.~65]{Z} випливає, що для функції $f_{\alpha_{0}}(t)\!:=\!(\varphi_{\alpha_{0}}\ast\Psi_{\beta})(t)$ має місце рівність
\begin{equation}\label{ochf1}
f_{\alpha_{0}}(t)\!=\!
\alpha_{0}\bigg(\!\sum\limits_{k=1}^{n}\!\psi(k)\cos\Big(kt-
\frac{\beta\pi}{2}\Big)\!+2\!\!\sum\limits_{k=n+1}^{2n-1}\!\psi(k)\Big(1-\frac{k}{2n}\Big)\cos\Big(kt-
\frac{\beta\pi}{2}\Big)\!\!\bigg).
\end{equation}

Розглянемо інтеграл \ \
$$
I:=\int\limits_{-\pi}^{\pi}\Big(f_{\alpha_{0}}(t)-Z_{n-1}^{s}(f_{\alpha_{0}};t)\Big)\sum\limits_{k=1}^{n-1}
\frac{\big(\psi(k)k^{s+\frac{1}{q'}}\big)^{q-1}}{k^{\frac{1}{q}}}\cos\Big(kt-\frac{\beta\pi}{2}\Big)dt.
$$
Враховуючи  (\ref{ochf1}), (\ref{sz}) та відому формулу
$$
\int\limits_{-\pi}^{\pi}\cos\Big(kt-\frac{\beta\pi}{2}\Big)\cos\Big(mt-\frac{\beta\pi}{2}\Big)dt={\left\{ {\begin{array}{l l}
0, & k\neq m, \ \ k,m\in \mathbb{N}, \ \ \beta\in \mathbb{R}, \\
 \pi, & k=m, \ \  k,m\in \mathbb{N}, \ \ \beta\in \mathbb{R},
\end{array}} \right.}
$$
інтеграл $I$ можна записати у вигляді \ \
$
I=
\frac{\alpha_{0}\pi}{n^{s}}\sum\limits_{k=1}^{n-1}\frac{\big(\psi(k)k^{s+\frac{1}{q'}}\big)^{q}}
{k}$, звідки одержуємо
\begin{equation}\label{bilche}
I\cdot\bigg(\sum\limits_{k=1}^{n-1}\frac{\big(\psi(k)k^{s+\frac{1}{q'}}\big)^{q}}
{k}\bigg)^{-\frac{1}{q'}}\! =\frac{\alpha_{0}\pi}{n^{s}}\bigg(\sum\limits_{k=1}^{n-1}\frac{\big(\psi(k)k^{s+\frac{1}{q'}}\big)^{q}}
{k}\bigg)^{\frac{1}{q}}\geq\frac{K_{\psi,s,q}^{(1)}}{n^{s}}\bigg(\int\limits_{1}^{n}\frac{\big(g_{s,q'}(t)\big)^{q}}{t}dt\bigg)^{\frac{1}{q}},
\end{equation}
тут і далі $ K^{(i)}_{\psi,s,q}$  --- додатні величини, що можуть залежати тільки від $\psi$, $s$ та $q$.

З іншого боку, для оцінки зверху величини $I$ застосуємо  нерівність Гельдера та проведемо міркування, використані при доведенні формули (42) з роботи \cite{SG}, внаслідок чого отримуємо, що при $1<q<\infty$
\begin{equation}\label{menche}
I\leq\Big\|f_{\alpha_{0}}(t)-Z_{n-1}^{s}(f_{\alpha_{0}};t)
\Big\|_{q}\bigg\|\sum\limits_{k=1}^{n-1}\frac{(g_{s,q'}(k))^{q-1}}
{k^{\frac{1}{q}}}\cos\Big(kt-\frac{\beta\pi}{2}\Big)\bigg\|_{q'}\leq
$$
$$
\leq K^{(2)}_{\psi,s,q}\Big\|f_{\alpha_{0}}(t)-Z_{n-1}^{s}(f_{\alpha_{0}};t)
\Big\|_{q}\Big(\sum\limits_{k=1}^{n-1}\frac{(g_{s,q'}(k))^{q}}{k}\Big)^{\frac{1}{q'}}, \ \frac{1}{q}+\frac{1}{q'}=1.
\end{equation}
Об'єднуючи (\ref{bilche})  та (\ref{menche}) будемо мати
\begin{equation}\label{drstipomn}
\frac{K^{(1)}_{\psi,s,q}}{n^{s}}\bigg(\int\limits_{1}^{n}\frac{\big(g_{s,q'}(t)\big)^{q}}{t}dt\bigg)^{\frac{1}{q}}\leq
I\cdot\bigg(\sum\limits_{k=1}^{n-1}\frac{\big(g_{s,q'}(k)\big)^{q}}
{k}\bigg)^{-\frac{1}{q'}} \leq
K^{(2)}_{\psi,s,q}\Big\|f_{\alpha_{0}}(t)-Z_{n-1}^{s}(f_{\alpha_{0}};t)
\Big\|_{q}.
\end{equation}
З  (\ref{drstipomn}) отримуємо
$$
{\cal E}\big(L^{\psi}_{\beta,1};
Z_{n-1}^{s}\big)_{q}\geq\Big\|f_{\alpha_{0}}(t)-Z_{n-1}^{s}(f_{\alpha_{0}};t)
\Big\|_{q}\geq \frac{K^{(3)}_{\psi,s,q}}{n^{s}}\bigg(\int\limits_{1}^{n}\frac{\big(g_{s,q'}(t)\big)^{q}}{t}dt\bigg)^{\frac{1}{q}}.
$$
Оцінку знизу величини ${\cal E}\big(L^{\psi}_{\beta,1};
Z_{n-1}^{s}\big)_{q}$ у випадку $g_{s,q'}\in \mathcal{Z}$, $1<q<\infty$ доведено.

 Як випливає  з теореми (2.2.1) роботи
\cite[с.~92]{S1},  метод $Z_{n-1}^{s}$ насичений з порядком
насичення ${n^{-s}}$,  тобто величини ${\cal E}\big(L^{\psi}_{\beta,1};
Z_{n-1}^{s}\big)_{q}$ не можуть прямувати до нуля швидше за порядком, ніж ${n^{-s}}$, звідки випливає,
 що для довільного $\beta\in \mathbb{R}$, $1<q<\infty$, за умови
$g_{s,q'}\in A^{-}$, виконується порядкова рівність (\ref{tp12p}). Теорему 2 доведено.

\textbf{Наслідок 1.} \emph{Нехай $1< q<\infty$,  $n\in \mathbb{N}$,
$g_{1,q'}(t):=\psi(t)t^{1+\frac{1}{q'}}$,   $\frac{1}{q}+\frac{1}{q'}=1$ і $\beta\in \mathbb{R}$. Тоді}

\emph{1. Якщо $\psi\in \Theta_{q'}$, $g_{1,q'}\in A^{+}$ і функція $1/\psi(t)$ опукла вгору або донизу на $[1,\infty)$, то}
$$
{\cal E}\big(L^{\psi}_{\beta,1};
\sigma_{n-1}\big)_{q}\asymp E_{n}\big(L^{\psi}_{\beta,1}
\big)_{q}\asymp\psi(n)n^{1-\frac{1}{q}}.
$$

\emph{2. Якщо   $g_{1,q'}\in \mathcal{Z}$,
   то}
$$
{\cal E}\big(L^{\psi}_{\beta,1};
\sigma_{n-1}\big)_{q}\asymp\frac{1}{n}
\Big(\int\limits_{1}^{n}\frac{\big(g_{1,q'}(t)\big)^{q}}{t}dt\Big)^{\frac{1}{q}}.
$$

\emph{3. Якщо  $g_{1,q'}\in A^{-}$, то}
$$
{\cal E}\big(L^{\psi}_{\beta,1};
\sigma_{n-1}\big)_{q}\asymp n^{-1}.
$$

\textbf{Наслідок 2.} \emph{Нехай  $r>1-\frac{1}{q}$, $1<q<\infty$, $s>0$, $\beta\in \mathbb{R}$, $n\in\mathbb{N}\backslash\{1\}$. Тоді}
\begin{equation}\label{vnad}
{\cal E}\left(W^{r}_{\beta,1};
Z_{n-1}^{s}\right)_{q}\asymp{\left\{ {\begin{array}{l l}
n^{-(r-1+\frac{1}{q})}, & 1-\frac{1}{q}<r<s+1-\frac{1}{q}; \\
n^{-s}\ln^{\frac{1}{q}}n, &r=s+1-\frac{1}{q};  \\
n^{-s}, &r>s+1-\frac{1}{q}.
\end{array}} \right.}
\end{equation}

 При $s=1$ із (\ref{vnad}) випливають результати одержані А.І. Камзоловим \cite{Kam}. Із (\ref{vnad}) і теореми 3.6 роботи  \cite[с.~47] {tem} випливає, що при $1-\frac{1}{q}<r<s+1-\frac{1}{q}$ і $\beta\in \mathbb{R}$
$$
{\cal E}\left(W^{r}_{\beta,1};
Z_{n-1}^{s}\right)_{q}\asymp E_{n}\left(W^{r}_{\beta,1}
\right)_{q}\asymp n^{-(r-1+\frac{1}{q})}, \ \ 1<q<\infty.
$$



\renewcommand{\refname}{}
\makeatletter\renewcommand{\@biblabel}[1]{#1.}\makeatother


\end{document}